\documentclass[10pt]{article}
\begin{document}
\title{ {\bf  On Potentially $K_6-C_5$ graphic
Sequences}
\thanks{   Project Supported by  NSF of Fujian(2008J0209 ),
Science and Technology Project of Fujian, Fujian Provincial Training
Foundation for "Bai-Quan-Wan Talents Engineering" , Project of
Fujian Education Department and Project of Zhangzhou Teachers
College(SK07009).}}
\author{{Zhenghua Xu, Chunhui Lai}\\
{\small Department of Mathematics, Zhangzhou Teachers College,}
\\{\small Zhangzhou, Fujian 363000,
 P. R. of CHINA.}\\{\small  xzh2002528@tom.com( Zhenghua Xu )}
 \\{\small   zjlaichu@public.zzptt.fj.cn(Chunhui Lai, Corresponding author)}
}

\date{}
\maketitle
\begin{center}
\begin{minipage}{4.1in}
\vskip 0.1in
\begin{center}{\bf Abstract}\end{center}
 \ \ { For given a graph $H$, a graphic sequence $\pi=(d_1,d_2,\cdots,d_n)$ is said to be potentially
 $H$-graphic if there exists a realization of $\pi$ containing $H$ as a subgraph. In this paper,
we characterize the potentially  $K_6-C_5$ -graphic sequences .}
\par
\end{minipage}
\end{center}
 \par
\def\CTeXPreproc{Created by ctex v0.2.12, don't edit!}
\section{Introduction}
\par
\baselineskip 14pt

\ \ \ \ Let $G$ be a simple undirected graph, without loops or
multiple edges. Any undefined notation follows that of Bondy and
Murty $[1]$. The set of all non-increasing nonnegative integer
sequence $\pi=(d_1,d_2,\cdots,d_n)$ is denoted by $NS_n$. A
sequence $\pi$ is called graphic if there is a (simple) graph $G$
of order n having degree sequence $\pi$ and such a graph $G$ is
referred as a realization of $\pi$. The set of all graphic
sequence in $NS_n$ is denoted by $GS_n$. A graphic sequence $\pi$
is potentially $H$-graphic if there is a realization of $\pi$
containing $H$ as a subgraph. Let $C_k$ and $P_k$ denote a cycle
on $k$ vertices and a path on $k+1$ vertices, respectively. Let
$\sigma(\pi)$ be the sum of all the terms of $\pi$.  Let $G-H$
denote the graph obtained from $G$ by removing the edges set
$E(H)$ where $H$ is a subgraph of $G$.  In the degree sequence,
$r^t$ means that there are $t$ vertices of degree $r$ in the
realization of the sequence.\par

  Tur\'{a}n number is the maximum number of edges of a graph
with $n$ vertices not containing $H$ as a subgraph and it is
denoted by $ex(n,H)$. In terms of graphic sequences, the number
$2ex(n,H)+2$ is the minimum even integer $l$ such that every
$n$-term graphical sequence $\pi$ with $\sigma (\pi)\geq l $ is
forcibly $H$-graphical. Gould, Jacobson and Lehel [3] considered
the following variation of the classical Tur\'{a}n-type extremal
problems: determine the smallest even integer $\sigma(H,n)$ such
that every n-term positive graphic sequence
$\pi=(d_1,d_2,\cdots,d_n)$ with $\sigma(\pi)\geq \sigma(H,n)$ has
a realization $G$ containing $H$ as a subgraph. They gave that
$\sigma(pK_2, n)=(p-1)(2n-p)+2$ for $p\ge 2$; $\sigma(C_4,
n)=2[{{3n-1}\over 2}]$ for $n\ge 4$.   Erd\"os,\ Jacobson and
Lehel [2] showed that $\sigma(K_k, n)\ge (k-2)(2n-k+1)+2$ and
conjectured  the equality holding and proved the conjecture is
true for $k=3$ and $n\geq 6$. The conjecture is confirmed in [3]
and [15]-[18]. Lai[9] determined $\sigma (K_4-e, n)$ for $n\geq
4$. Yin[25] and Lai[12] determined $\sigma(K_{1,1,3},n)$ for
$n\geq 5$ independently.
 Yin,Li and Mao[27] determined $\sigma(K_{r+1}-e,n)$ for $r\geq 3,$ $r+1\leq n
\leq 2r$ and $\sigma(K_5-e,n)$ for $n\geq5$. Yin and Li[26] gave a
good method of determining  $\sigma(K_{r+1}-e,n)$ for $r\geq2$ and
$n\geq3r^2-r-1$ (In fact, Yin and Li[26] also determined
$\sigma(K_{r+1}-ke,n)$ for $r\geq2$ and $n\geq3r^2-r-1$). After
reading[26], Yin [31] determined $\sigma (K_{r+1}-K_{3}, n)$ for
    $n\geq 3r+5, r\geq 3$.Lai [10-11] determined
    $\sigma (K_{5}-C_{4}, n),\sigma (K_{5}-P_{3}, n)$ and
    $\sigma (K_{5}-P_{4}, n)$ where  $n\geq 5$. Determining $\sigma(K_{r+1}-H,n)$ where $H$
    is a tree on 4 vertices ,is more useful than a cycle on 4
    vertices (for example, although $C_4 \not\subset C_i$,  $P_3 \subset C_i$ for $i\geq 5$).
    After reading[26] and [31], Lai and Hu[13]
    determined  $\sigma (K_{r+1}-H, n)$ for
    $n\geq 4r+10, r\geq 3, r+1 \geq k \geq 4$ where $H$ is a graph on $k$
    vertices which
    contains a tree on $4$ vertices but
     not contain a cycle on $3$ vertices and $\sigma (K_{r+1}-P_2, n)$ for
    $n\geq 4r+8, r\geq 3$.  Recently, Lai $[14]$ also determined $\sigma(K_{r+1}-Z_4,n)$,
$\sigma(K_{r+1}-(K_4-e),n)$, $\sigma(K_{r+1}-K_4,n)$ for $n\geq
5r+16$, $r\geq 4$ and $\sigma(K_{r+1}-Z,n)$ for $n\geq5r+19$,
$r+1\geq k\geq5$, $j\geq 5$ where $Z$ is a graph on $k$ vertices and
$j$ edges which contains a graph $Z_4$ but not contain a cycle on 4
vertices.
\par
  A harder question is to characterize the potentially
 $H$-graphic sequences without zero terms.  Luo [20] characterized the potentially
 $C_k$-graphic sequences for each $k=3,4,5$ and  Luo and Warner [21] characterized the potentially
 $K_4$-graphic sequences. Lately, Eschen and Niu [22] characterized the potentially
 $K_4-e$-graphic sequences.  Yin and Chen [29] characterized the
 potentially $K_{r,s}$-graphic sequences for $r=2,s=3$ and
 $r=2,s=4$. Yin et al. [32] characterized the
 potentially $K_5-e$ and $K_6$-graphic sequences. Recently, Hu et al. [4-6] characterized the potentially
 $K_5-C_4$, $K_5-P_4$ and $K_5-E_3$-graphic
 sequences where $E_3$ denotes graphs with 5 vertices and 3 edges and the potentially
 $K_{3,3}$, $K_6-C_6$-graphic sequences in [7]. Yin [33] characterized the potentially
$K_6-K_3$-graphic sequences. In this paper, we have characterized
the potentially  $K_6-C_5$ -graphic sequences.
\par
\section{Preparations}\par

  \ \ \ \  Let $\pi=(d_1,\cdots,d_n)\epsilon NS_n,1\leq k\leq n$. Let
    $$ \pi_k^{\prime\prime}=\left\{
    \begin{array}{ll}(d_1-1,\cdots,d_{k-1}-1,d_{k+1}-1,
    \cdots,d_{d_k+1}-1,d_{d_k+2},\cdots,d_n), \\ \mbox{ if $d_k\geq k,$}\\
    (d_1-1,\cdots,d_{d_k}-1,d_{d_k+1},\cdots,d_{k-1},d_{k+1},\cdots,d_n),
     \\ \mbox{if $d_k < k.$} \end{array} \right. $$
  Denote
  $\pi_k^\prime=(d_1^\prime,d_2^\prime,\cdots,d_{n-1}^\prime)$, where
  $d_1^\prime\geq d_2^\prime\geq\cdots\geq d_{n-1}^\prime$ is a
  rearrangement of the $n-1$ terms of $\pi_k^{\prime\prime}$. Then
  $\pi_k^{\prime}$ is called the residual sequence obtained by
  laying off $d_k$ from $\pi$. For simplicity, we denote $\pi_n^\prime$ by $\pi^\prime$ in this paper.
 \par
  If $S\subseteq NS_n$, then $\pi\in S$ if and only if $\pi$ is an
  element of $S$ by rearranging the terms of $\pi$. \par

  Let $\pi_1=(d_1^1,d_2^1,\cdots,d_n^1)$$\in NS_n$, $\pi_2=(d_1^2,d_2^2,\cdots,d_m^2)$$\in
  NS_n$
  and $n\geq m$, $d_i^1\geq d_i^2$$(i=1,2,\cdots,m)$. Let
  $\pi^*=$$\pi_1-\pi_2=$$(d_1^1-d_1^2,d_2^1-d_2^2,\cdots,d_m^1-d_m^2,d^1_{m+1},\cdots,d^1_n)$.\par

   For a non-increasing positive integer sequence $\pi=(d_1,d_2,\cdots,d_n)$,
   we write $m(\pi)$ and $h(\pi)$ to denote the largest
positive terms of $\pi$ and the smallest positive terms of $\pi$,
respectively. We need the following the results.
\par
    {\bf Theorem 2.1 [3]}\ \   If $\pi=(d_1,d_2,\cdots,d_n)$ is a graphic
 sequence with a realization $G$ containing $H$ as a subgraph,
 then there exists a realization $G^\prime$ of $\pi$ containing $H$ as a
 subgraph so that the vertices of $H$ have the largest degrees of
 $\pi$.\par
 \par
    {\bf Theorem 2.2 [19]}\ \   If $\pi=(d_1,d_2,\cdots,d_n)$ is a
 sequence of nonnegative integers where $1\leq m(\pi)\leq2$,
 $h(\pi)=1$ and $\sigma(\pi)$ is  even, then $\pi$ is graphic.
\par
    {\bf Lemma 2.3 (Kleitman and Wang [8])}\ \   $\pi$ is
graphic if and only if $\pi^\prime$ is graphic.
 \par
    The following corollary is obvious.\par
\par
    {\bf Corollary 2.4}\ \    Let $H$ be a simple graph. If $\pi^\prime$ is
 potentially $H$-graphic, then $\pi$ is
 potentially $H$-graphic.\par
    {\bf Lemma 2.5 [24]}\ \  Let $G$ be a simple undirected graph with the vertex $v_i$ and
    $d_G(v^\prime_i)=d^\prime_i$($i=1,2,\cdots,m$), $d_1^\prime \geq d_2^\prime \geq  \cdots \geq  d_m^\prime$.
    If $\pi=(d_1,d_2,\cdots,d_n)$ $\in
    NS_n$ is potentially $G$-graphic and $n-m \geq  (d_1-d^\prime_m)$, then there exists a realization $H$ of
$\pi$ containing $G$ as a
 subgraph so that the vertices of $G$ have the largest degrees of
 $\pi$ in proper order, that is, if the vertices of $G$ are $v^\prime_i$, $i=1,2,\cdots,m$$(m\leq
 n)$,
 $d_G(v^\prime_i)=d^\prime_i$ and $d_1^\prime \geq d_2^\prime \geq  \cdots\geq  d_m^\prime$, then
 there exists a realization $H$ of $\pi$ containing $G$ as a
 subgraph so that the vertices of $H$ are $v_i$, $i=1,2,\cdots,n$,
 $d_H(v_i)=d_i$, $G\subseteq H[v_1,v_2,\cdots,v_m]$ and the vertices $v_i^\prime$ are
  corresponding to the vertices $v_i$, $i=1,2,\cdots,m$.\par

    {\bf Lemma 2.6 [23]}\ \ Let
    $\pi=(4^{x_1},3^{x_2},2^{x_3},1^{x_4})$ where
    $\sigma(\pi)$ is even, $x_1+x_2+x_3+x_4=n$ and $n\geq 1$. Then $\pi\in GS_n$
    if and only if $\pi \notin S$, where
    $S=$
    $\{$$(2)$,$(2^2)$,$(3,1)$,$(3^2)$,$(3,2,1)$,$(3^2,2)$,\,$(3^3,1)$,\,$(3^2,1^2)$,\,$(4)$,\,\,$(4,1^2)$,$(4,2)$,
    $(4,2^2)$,\\$(4,2^3)$,$(4,2,1^2)$,$(4,3^2)$,$(4,3^2,2)$,$(4,3,1)$,$(4,3,1^3)$,\,$(4,3^2,1^2)$,\,$(4,3,2,1)$,
    $(4^2)$,\\$(4^2,1^2)$,\,\,\,$(4^2,1^4)$,\,\,\,$(4^2,2,1^2)$,\,\,\,$(4^2,2)$,\,\,\,$(4^2,2^2)$,\,\,\,$(4^2,3^2)$,\,\,\,$(4^2,3,1)$,\,\,\,$(4^2,3,1^3)$,\\
    $(4^2,3,2,1)$,\,\,$(4^3)$,\,\,$(4^3,1^2)$,\,\,$(4^3,2,1^2)$,\,\,$(4^3,1^4)$,\,\,$(4^3,2)$,\,\,$(4^3,2^2)$,\,\,$(4^3,3,1)$,$(4^4)$,\\
    $(4^4,1^2)$,$(4^4,2)$$\}$.\par

    In order to prove Theorem 3.1, at first, we prove Lemma 2.7 to Lemma
    2.11. Let $\pi_1^1=(5,3^5)$ be a graphic sequence of $K_6-C_5$ and $\pi^*$$=\pi-\pi_1^1$.\par

    {\bf Lemma 2.7 }\ \ $(1)$ If $\pi=(5^4,4^{n-4})$ where  $n\geq 6$ and $\sigma(\pi)$ is even,
    then $\pi$ is potentially $K_6-C_5$-graphic;\par
     $(2)$ If $\pi=(5^3,4^i,3^j,2^{n-3-i-j})$ where $n\geq 6$, $i+j\geq 3$ and $\sigma(\pi)$ is even, then $\pi$
  is potentially $K_6-C_5$-graphic if and only if $\pi\neq (5^3,3^3)$.\par

   {\bf Proof:}  $(1)$ If $\pi=(5^4,4^{n-4})$ where $n\geq 6$, then
   $\pi^*=(2,2,2,1,1,4^{n-6})$. Let $\pi_1=(4^{n-6})$.
    If $\pi_1=(4^{n-6})$$\notin S$, according to
 Lemma 2.6, then $\pi_1$ is graphic. Let $G$ be a realization of
 $\pi_1=(4^{n-6})$, then $(K_6-e)$$\bigcup G$ is a realization of
 $\pi=(5^4,4^{n-4})$, hence  $\pi$ is potentially $K_6-C_5$-graphic.
 If $\pi_1=(4^{n-6})$$\in S$, then
 $\pi_1=$$(4)$, $(4^2)$, $(4^3)$, $(4^4)$.
 Therefore, $\pi=$$(5^4,4^3)$, $(5^4,4^4)$, $(5^4,4^5)$, $(5^4,4^6)$.
 It is easy to check that all of these are potentially $K_6-C_5$-graphic.\par
    $(2)$ If $\pi=(5^3,4^i,3^j,2^{n-3-i-j})$ where $n\geq 6$ and $i+j\geq 3$, then
 we consider the following cases.

\textbf{Case 1:} $i=0$ \par
   If $\pi=(5^3,3^j,2^{n-3-j})$$(j\geq 3)$, then $\pi^*=(2^2,3^{j-3},2^{n-3-j})$.

   If j=3, then $\pi^*=(2^2,2^{n-6})$. If $n=6$, then $\pi=(5^3,3^3)$.
   It is easy to check that $\pi=$$(5^3,3^3)$
 is not potentially $K_6-C_5$-graphic.
  If $n\geq 7$, then let $\pi_1=(1,1,2^{n-7})$.
  Since $\pi_1=(1,1,2^{n-7})$$\notin S$, according to
 Lemma 2.6, then $\pi_1$  is graphic. Firstly, let $G_1$ be a realization of
 $\pi_1=(1,1,2^{n-7})$ and $u_1$ and $u_2$ denote two vertices of degree 1 in $G_1$.
 Secondly, we add a vertex $u_3$ to $G_1$. Let $u_3$ be
 adjacent to $u_1$ and $u_2$. Then we obtain a graph $G_2$ and
 $G_2$ is a realization of $\pi^*=(2^2,2^{n-6})$($n\geq 7$).
 Finally, let $v_1$ and $v_2$ denote two vertices of degree 3 and be not adjacent  in  $K_6-C_5$.
 We replace $v_1$ and $u_1$ with a single vertex whose incident
 edges are the edges that were incident to $v_1$ or $u_1$
 and we also replace $v_2$ and $u_2$ with a single vertex whose incident
 edges are the edges that were incident to $v_2$ or $u_2$ in $G_2$$\bigcup$($K_6-C_5$).
 Then  we obtain a graph $G$ and $G$ is a realization of
 $\pi=(5^3,3^3,2^{n-6})$, hence $\pi$ is potentially
 $K_6-C_5$-graphic. We similarly do as follows. \par

  If $j\geq 4$, then let $\pi_1=(1,2,3^{j-4},2^{n-3-j})$.
  If $\pi_1=(1,2,3^{j-4},2^{n-3-j})$$\notin S$, according to
 Lemma 2.6, then $\pi_1$  is graphic, hence  $\pi$ is potentially
 $K_6-C_5$-graphic.
 If $\pi_1=(1,2,3^{j-4},2^{n-3-j})$$\in S$, then
 $\pi_1=$$(3,2,1)$.
 Therefore, $\pi=$$(5^3,3^5)$.
 It is easy to check that $\pi=$$(5^3,3^5)$
 is potentially $K_6-C_5$-graphic .\par

 \textbf{Case 2:} $i=1$ \par
   If $\pi=(5^3,4,3^j,2^{n-4-j})$$(j\geq 2)$, then
   $\pi^*=(2,2,1,3^{j-2},2^{n-4-j})$. Let
   $\pi_1=(1,3^{j-2},2^{n-4-j})$. If $\pi_1=(1,3^{j-2},2^{n-4-j})$$\notin S$, according to
 Lemma 2.6, then $\pi_1$  is graphic, hence  $\pi$ is potentially
 $K_6-C_5$-graphic.
 If $\pi_1=(1,3^{j-2},2^{n-4-j})$$\in S$, then
 $\pi_1=$$(3,1)$, $(3,2,1)$, $(3^3,1)$. Therefore, $\pi=$$(5^3,4,\\3^3)$,$(5^3,4,3^3,2)$,$(5^3,4,3^5)$.
 It is easy to check that all of these are potentially $K_6-C_5$-graphic .\par

\textbf{Case 3:} $i=2$ \par
   If $\pi=(5^3,4^2,3^j,2^{n-5-j})$$(j\geq 1)$, then
$\pi^*=(2,2,1,1,3^{j-1},2^{n-5-j})$. Let
$\pi_1=(3^{j-1},2^{n-5-j})$. If $\pi_1=(3^{j-1},2^{n-5-j})$$\notin
S$, according to Lemma 2.6, then $\pi_1$  is graphic, hence  $\pi$
is potentially $K_6-C_5$-graphic. If
$\pi_1=(3^{j-1},2^{n-5-j})$$\in S$, then $\pi_1=$$(2)$, $(2^2)$,
$(3^2)$, $(3^2,2)$. Therefore, $\pi=$ $(5^3,4^2,\\3,2)$,
$(5^3,4^2,3,2^2)$, $(5^3,4^2,3^3)$, $(5^3,4^2,3^3,2)$. It is easy
to check that all of these are potentially $K_6-C_5$-graphic .\par

\textbf{Case 4:} $i\geq 3$ \par
   If \ $\pi=(5^3,4^i,3^j,2^{n-3-i-j})$ $(j\geq0 )$, \ then \ $\pi^*=(2,2,1,1,1,4^{i-3},3^j,\\2^{n-3-i-j})$. Let
   $\pi_1=(1,4^{i-3},3^j,2^{n-3-i-j})$.
   If $\pi_1=(1,4^{i-3},3^j,2^{n-3-i-j})$ $\notin S$,
   according to
 Lemma 2.6, then $\pi_1$  is graphic,
hence  $\pi$ is potentially
 $K_6-C_5$-graphic.
 If $\pi_1=(1,4^{i-3},3^j,2^{n-3-i-j})$$\in S$, then
 $\pi_1=$$(3,1)$, $(3,2,1)$, $(3^3,1)$,
 $(4,3,1)$, $(4,3,2,1)$, $(4^2,3,1)$, $(4^2,3,2,1)$, $(4^3,3,1)$.
 Therefore, $\pi=$$(5^3,4^3,3)$,
 $(5^3,4^3,3,2)$,$(5^3,4^3,3^3)$,$(5^3,4^4,3)$,$(5^3,4^4,3,2)$,
 $(5^3,4^5,3)$,\\$(5^3,4^5,3,2)$,
 $(5^3,4^6,3)$. It is easy to check that all of these are potentially
 $K_6-C_5$.\par

    {\bf Lemma 2.8 }\ \ If $\pi=(5^2,4^i,3^j,2^{n-2-i-j})$
where   $i+j\geq 4$, $n\geq 6$ and $\sigma(\pi)$ is even, then
$\pi$ is potentially $K_6-C_5$-graphic if and only if
$\pi\neq$$(5^2,3^4)$,\\$(5^2,3^4,2)$,$(5^2,3^6)$,$(5^2,4,3^4)$.\par

    {\bf Proof:} If $\pi=(5^2,4^i,3^j,2^{n-2-i-j})$ where $i+j\geq 4$ and
$n\geq 6$, then  we consider the following cases. \par

\textbf{Case 1:} $i=0$ \par
 If $\pi=(5^2,3^j,2^{n-2-j})$$(j\geq 4)$, then $\pi^*=(2,3^{j-4},2^{n-2-j})$.
 If $\pi^*\notin S$, according to
Lemma 2.6, then $\pi^*$ is graphic, hence $\pi$ is potentially
$K_6-C_5$-graphic. If $\pi^*\in S$, then
$\pi^*=$$(2)$,$(2^2)$,$(3^2,2)$. Therefore,
$\pi=(5^2,3^4)$,$(5^2,3^4,2)$,$(5^2,3^6)$. It is easy to check
that these are not potentially $K_6-C_5$-graphic .\par

\textbf{Case 2:} $i=1$ \par
 If $\pi=(5^2,4,3^j,2^{n-3-j})$$(j\geq 3)$, then
 $\pi^*=(2,1,3^{j-3},2^{n-3-j})$. Let
 $\pi_1=(1,3^{j-3},2^{n-3-j})$.
 If $\pi_1=(1,3^{j-3},2^{n-3-j})$$\notin S$, according to
 Lemma 2.6,then $\pi_1$ is graphic, hence  $\pi$ is potentially $K_6-C_5$-graphic.
 If $\pi_1=(1,3^{j-3},2^{n-3-j})$$\in S$,then
 $\pi_1=$$(3,1)$,$(3,2,1)$,$(3^3,1)$.Therefore, $\pi=$$(5^2,4,3^4)$,
 $(5^2,4,3^4,2)$,$(5^2,4,3^6)$.
 It is easy to check that others
 are  potentially $K_6-C_5$-graphic except for $(5^2,4,3^4)$.\par

\textbf{Case 3:} $i=2$ \par
 If $\pi=(5^2,4^2,3^j,2^{n-4-j})$$(j\geq 2)$, then
 $\pi^*=(2,1,1,3^{j-2},2^{n-4-j})$. Let
 $\pi_1=(3^{j-2},2^{n-4-j})$.
 If $\pi_1=(3^{j-2},2^{n-4-j})$$\notin S$, according to
 Lemma 2.6, then $\pi_1$  is graphic, hence $\pi$ is
potentially $K_6-C_5$-graphic.
 If $\pi_1=(3^{j-2},2^{n-4-j})$$\in S$,then
 $\pi_1=$$(2)$,$(2^2)$,$(3^2)$,$(3^2,2)$.Therefore,$\pi=(5^2,4^2,3^2,2)$,
 $(5^2,4^2,3^2,2^2)$, $(5^2,4^2,3^4)$, $(5^2,4^2,3^4,2)$. It is easy to check that
 all of these are potentially $K_6-C_5$-graphic .\par

\textbf{Case 4:} $i=3$ \par
 If $\pi=(5^2,4^3,3^j,2^{n-5-j})$$(j\geq 1)$, then
 $\pi^*=(2,1,1,1,3^{j-1},2^{n-5-j})$. Let
 $\pi_1=(1,3^{j-1},2^{n-5-j})$.
 If $\pi_1=(1,3^{j-1},2^{n-5-j})$$\notin S$, according to
 Lemma 2.6, then $\pi_1$ is graphic, hence $\pi$ is
potentially
 $K_6-C_5$-graphic.
 If $\pi_1=(1,3^{j-1},2^{n-5-j})$$\in S$, then $\pi_1=$$(3,1)$, $(3,2,1)$, $(3^3,1)$.
 Therefore, $\pi=$$(5^2,4^3,3^2)$,
 $(5^2,4^3,3^2,2)$,$(5^2,4^3,3^4)$. It is easy to check that
 all of these are potentially $K_6-C_5$-graphic .\par

\textbf{Case 5:} $i\geq 4$ \par
 If \ $\pi$ $=$ $(5^2,4^i,3^j,2^{n-2-i-j})$ $(j\geq 0)$, then\ \
 $\pi^*=(2,1,1,1,1,4^{i-4},3^j,\\2^{n-2-i-j})$. Let
 $\pi_1=(4^{i-4},3^j,2^{n-2-i-j})$.
 If $\pi_1=(4^{i-4},3^j,2^{n-2-i-j})$$\notin S$, according to
 Lemma 2.6, then $\pi_1$  is graphic, hence  $\pi$ is
potentially $K_6-C_5$-graphic.
 If $\pi_1=(4^{i-4},3^j,2^{n-2-i-j})$$\in S$, then
 $\pi_1=$$(2)$, $(2^2)$, $(3^2)$, $(3^2,2)$, $(4)$,
 $(4,2)$, $(4,2^2)$, $(4,2^3)$, $(4,3^2)$,
 $(4,3^2,2)$, $(4^2)$, $(4^2,2)$, $(4^2,2^2)$,
 $(4^2,3^2)$, $(4^3)$, $(4^3,2)$, $(4^3,2^2)$,
 $(4^4)$, $(4^4,2)$.
 Therefore, $\pi=$$(5^2,4^4,2)$, $(5^2,4^4,2^2)$, $(5^2,4^4,3^2)$,
 $(5^2,4^4,3^2,2)$, $(5^2,4^5)$, $(5^2,4^5,2)$, $(5^2,4^5,2^2)$, $(5^2,4^5,2^3)$, $(5^2,4^5,\\3^2)$,
 $(5^2,4^5,3^2,2)$, $(5^2,4^6)$, $(5^2,4^6,2)$, $(5^2,4^6,2^2)$, $(5^2,4^6,3^2)$,
 $(5^2,4^7)$, $(5^2,4^7,\\2)$, $(5^2,4^7,2^2)$, $(5^2,4^8)$, $(5^2,4^8,2)$.
 It is easy to check that
 all of these are potentially $K_6-C_5$-graphic .\par

 {\bf Lemma 2.9 }\ \ If $\pi=(5,4^i,3^j,2^k,1^{n-i-j-k-1})$ where $i+j\geq 5$, $n\geq 6$
 and $\sigma(\pi)$ is even, then $\pi$ is potentially
$K_6-C_5$-graphic if and only if $\pi\neq$$(5,3^5,2)$,
$(5,3^5,2^2)$, $(5,3^7)$, $(5,3^6,1)$, $(5,3^6,2,1)$, $(5,3^7,2)$,
$(5,3^8,1)$, $(5,3^7,1^2)$, $(5,4,3^5)$, $(5,4,3^5,2)$,
$(5,4,3^7)$, $(5,4,3^6,1)$, $(5,4^2,3^5)$.\par

{\bf Proof:} If $\pi=(5,4^i,3^j,2^k,1^{n-i-j-k-1})$ where $i+j\geq
5$ and $n\geq 6$, then  we consider the following cases. \par

\textbf{Case 1:} $i=0$\par
 If $\pi=(5,3^j,2^k,1^{n-j-k-1})$$(j\geq 5)$, then
 $\pi^*=(3^{j-5},2^k,1^{n-j-k-1})$.
 If $\pi^*\notin S$, according to
 Lemma 2.6, then $\pi^*$  is graphic, hence
$\pi$ is potentially $K_6-C_5$-graphic.
 If $\pi^*\in S$, then
 $\pi^*=$$(2)$, $(2^2)$, $(3^2)$, $(3,1)$,
 $(3,2,1)$,$(3^2,2)$,$(3^3,1)$,$(3^2,1^2)$.
 Therefore, $\pi=$$(5,3^5,2)$, $(5,3^5,2^2)$,
 $(5,3^7)$, $(5,3^6,1)$, $(5,3^6,2,1)$, $(5,3^7,2)$, $(5,3^8,1)$,
 $(5,3^7,1^2)$.
 It is easy to check that
 these are not potentially $K_6-C_5$-graphic.\par

\textbf{Case 2:} $i=1$\par
 If $\pi=(5,4,3^j,2^k,1^{n-2-j-k})$$(j\geq 4)$, then
 $\pi^*=(1,3^{j-4},2^k,1^{n-2-j-k})$.
 If $\pi^*\notin S$, according to
 Lemma 2.6, then $\pi^*$ is graphic, hence
$\pi$ is potentially
 $K_6-C_5$-graphic.
 If $\pi^*=(1,3^{j-4},2^k,1^{n-2-j-k})$$\in S$, then
 $\pi^*=$$(3,1)$,$(3,2,1)$,\\$(3^3,1)$,$(3^2,1^2)$.
 Therefore,$\pi=$$(5,4,3^5)$,$(5,4,3^5,2)$,
 $(5,4,3^7)$,$(5,4,3^6,1)$. It is easy to check that
 these are not potentially $K_6-C_5$-graphic.\par

\textbf{Case 3:} $i=2$\par
 If $\pi=(5,4^2,3^j,2^k,1^{n-3-j-k})$$(j\geq 3)$,then $\pi^*=(1,1,3^{j-3},2^k,1^{n-3-j-k})$.
 Let $\pi_1=(3^{j-3},2^k,1^{n-3-j-k})$.
 If $\pi_1=(3^{j-3},2^k,1^{n-3-j-k})$$\notin S$, according to
 Lemma 2.6, then $\pi_1$  is graphic, hence
$\pi$ is potentially
 $K_6-C_5$-graphic. If $\pi_1=(3^{j-3},2^k,1^{n-3-j-k})$$\in S$, then
 $\pi_1=$$(2)$, $(2^2)$, $(3^2)$, $(3,1)$, $(3,2,1)$,
 $(3^2,2)$, $(3^3,1)$, $(3^2,1^2)$.
 Therefore,
 $\pi=$$(5,4^2,3^3,2)$, $(5,4^2,3^3,2^2)$,
  $(5,4^2,3^5)$, $(5,4^2,3^4,1)$, $(5,4^2,3^4,2,1)$,
 $(5,4^2,3^5,2)$,$(5,4^2,3^6,1)$,$(5,4^2,3^5,1^2)$.
 It is easy to check that others
 are potentially $K_6-C_5$-graphic except for $(5,4^2,3^5)$.\par

\textbf{Case 4:} $i=3$\par
 If \ $\pi=$ \ $(5,4^3,3^j,2^k,1^{n-4-j-k})$ \ $(j\geq 2)$, \ then \ $\pi^*=$ \ $(1,1,1,3^{j-2},2^k,\\1^{n-4-j-k})$.
Let $\pi_1=(1,3^{j-2},2^k,1^{n-4-j-k})$. If
$\pi_1=(1,3^{j-2},2^k,1^{n-4-j-k})$ $\notin S$, according to Lemma
2.6, then $\pi_1$  is graphic, hence $\pi$ is potentially
$K_6-C_5$-graphic. If $\pi_1=(1,3^{j-2},2^k,1^{n-4-j-k})$$\in S$,
then\,\,
 $\pi_1=$$(3,1)$, $(3,2,1)$, $(3^3,1)$,
 $(3^2,1^2)$.
 Therefore, $\pi=$$(5,4^3,3^3)$, $(5,4^3,3^3,2)$,
 $(5,4^3,3^5)$,$(5,4^3,3^4,1)$.
 It is easy to check that
 all of these are potentially $K_6-C_5$-graphic .\par

\textbf{Case 5:} $i=4$\par
 If \ $\pi=$ \ $(5,4^4,3^j,2^k,1^{n-5-j-k})$ \ $(j\geq 1)$, then
 $\pi^*=$ \ $(1,1,1,1,3^{j-1},2^k,\\1^{n-5-j-k})$. Let
 $\pi_1=(3^{j-1},2^k,1^{n-5-j-k})$.
If $\pi_1=(3^{j-1},2^k,1^{n-5-j-k})$$\notin S$,
 according to
 Lemma 2.6, then $\pi_1$  is graphic,
hence $\pi$ is potentially
 $K_6-C_5$-graphic.
 If $\pi_1=(3^{j-1},2^k,1^{n-5-j-k})$$\in S$, then
 $\pi_1=$$(2)$, $(2^2)$, $(3^2)$, $(3,1)$,
 $(3,2,1)$, $(3^2,2)$, $(3^3,1)$, $(3^2,1^2)$.
 Therefore, $\pi=$$(5,4^4,3,2)$, $(5,4^4,3,2^2)$,
 $(5,4^4,3^3)$, $(5,4^4,3^2,1)$, $(5,4^4,3^2,2,1)$, $(5,4^4,3^3,2)$,
 $(5,4^4,3^4,1)$, $(5,4^4,3^3,1^2)$.
 It is easy to check that
 all of these are potentially $K_6-C_5$-graphic .\par

\textbf{Case 6:} $i\geq 5$\par
 If $\pi=(5,4^i,3^j,2^k,2^{n-1-i-j-k})$($j\geq 0$ and $\sigma(\pi)$ is even),
 then $\pi$ is a graphic sequence,
 hence
 $\pi^*=(1,1,1,1,1,4^{i-5},3^j,2^k,1^{n-1-i-j-k})$.
 Let $\pi_1=(1,4^{i-5},3^j,2^k,1^{n-1-i-j-k})$.
 If \,\,\,\,$\pi_1=(1,4^{i-5},3^j,2^k,1^{n-1-i-j-k})$$\notin S$, according to
 Lemma 2.6, then $\pi_1$
 is graphic, hence $\pi$ is
potentially $K_6-C_5$-graphic.
 If $\pi_1=(1,4^{i-5},3^j,2^k,1^{n-1-i-j-k})$$\in
 S$, then $\pi_1=$$(3,1)$,\,\,  $(3,2,1)$, $(3^3,1)$,\,\, $(3^2,1^2)$,\,\, $(4,1^2)$,\,\,
$(4,2,1^2)$,\, $(4,3,1)$,\, $(4,3,1^3)$,\, $(4,3^2,1^2)$, \\
$(4,3,2,1)$,$(4^2,1^2)$,$(4^2,1^4)$,$(4^2,2,1^2)$,
 $(4^2,3,1)$,$(4^2,3,1^3)$,$(4^2,3,2,1)$,$(4^3,1^2)$,
 $(4^3,2,1^2)$, $(4^3,1^4)$, $(4^3,3,1)$, $(4^4,1^2)$.
 Therefore, $\pi=$$(5,4^5,3)$, $(5,4^5,3,2)$,
 $(5,4^5,3^3)$, $(5,4^5,3^2,1)$, $(5,4^6,1)$, $(5,4^6,2,1)$, $(5,4^6,3)$, $(5,4^6,3,1^2)$,
 $(5,4^6,\\3^2,1)$, $(5,4^6,3,2)$, $(5,4^7,1)$,
 $(5,4^7,1^3)$, $(5,4^7,2,1)$, $(5,4^7,3)$, $(5,4^7,3,1^2)$,
 $(5,4^7,3,2)$, $(5,4^8,1)$, $(5,4^8,2,1)$, $(5,4^8,1^3)$,
 $(5,4^8,3)$, $(5,4^9,1)$.
 It is easy to check that all of these are potentially $K_6-C_5$-graphic .\par

  {\bf Lemma 2.10}\ \ If $\pi=(d_1,d_2,3^{n-2})$ where $d_1\geq 5$,
 $d_2\geq 3$, $n\geq 6$ and $\sigma(\pi)$ is even, then \par
 $(1)$ $\pi\in GS_n$;\par
 $(2)$  $\pi$ is potentially $K_6-C_5$-graphic if  the following conditions hold:\par
  $(i)$ If \,$n\geq 7$ and $n$ is even,  then $\pi\neq ((n-1)^2,3^{n-2})$
  and \,$\pi\neq ((n-2)^2,3^{n-2})$;\par

  $(ii)$ If $n\geq 7$ and $n$ is odd, then $\pi\neq
  ((n-1),(n-2),3^{n-2})$;\par
  $(iii)$ $\pi\neq$$(5,3^7)$,$(5,4,3^5)$,$(5,4,3^7)$,$(5^2,3^4)$,$(5^2,3^6)$,
  $(6,3^6)$,$(6,3^8)$,$(7,3^7)$,\\$(6,4,3^6)$.

{\bf Proof:} $(1)$ is obvious.\par
  $(2)$ If $\pi=(d_1,d_2,3^{n-2})$ where $d_1\geq 5$,
 $d_2\geq 3$, $n\geq 6$ and $\sigma(\pi)$ is even,
 then $\pi$ is potentially $K_6-C_5$-graphic when
$\pi^*=(d_1-5,d_2-3,3^{n-6})$ is graphic  and two vertices of
degree $d_1-5$ and $d_2-3$ are not adjacent. We consider the
following cases.\par

\textbf{Case 1:}  $d_1-5=0$\par
  If $\pi=(5,d_2,3^{n-2})$, then $\pi^*=(d_2-3,3^{n-6})$.

  If $d_2-3=0$, then $\pi^*=(3^{n-6})$. If $\pi^*\notin S$, according to
 Lemma 2.6, then $\pi^*$ is graphic, hence $\pi$ is potentially
 $K_6-C_5$-graphic.
 If $\pi^*\in S$, then
 $\pi^*=$$(3^2)$. Therefore, $\pi=(5,3^7)$. It is easy to check that
 $\pi=(5,3^7)$ is not potentially $K_6-C_5$-graphic.\par

 If $d_2-3=1$, then $\pi^*=(1,3^{n-6})$. If $\pi^*\notin S$, according to
 Lemma 2.6, then $\pi^*$ is graphic, hence $\pi$ is potentially
 $K_6-C_5$-graphic.
 If $\pi^*\in S$,then $\pi^*=$$(3,1)$,$(3^3,1)$. Therefore, $\pi=(5,4,3^5)$,$(5,4,3^7)$.
It is easy to check that
 $\pi=(5,4,3^5)$ and $(5,4,3^7)$ are not potentially $K_6-C_5$-graphic.\par

If $d_2-3=2$, then $\pi^*=(2,3^{n-6})$. If $\pi^*\notin S$,
according to Lemma 2.6, then $\pi^*$ is graphic, hence $\pi$ is
potentially
 $K_6-C_5$-graphic.
 If $\pi^*\in S$, then
 $\pi^*=$$(2)$,$(3^2,2)$. Therefore, $\pi=(5^2,3^4)$,$(5^2,3^6)$. It is easy to check that
 $\pi=(5^2,3^4)$ and $(5^2,3^6)$ are not potentially $K_6-C_5$-graphic.\par

\textbf{Case 2:}  $d_1-5\neq 0$\par
  If $\pi=(d_1,d_2,3^{n-2})$, then $\pi^*=(d_1-5,d_2-3,3^{n-6})$.

  If $d_2-3=0$, then $\pi^*=(d_1-5,3^{n-6})$. Let
   $\pi_1=$$(3^{(n-6)-(d_1-5)},2^{d_1-5})$.
 If $\pi_1\notin S$, according to
 Lemma 2.6, then $\pi_1$ is graphic, hence  $\pi$ is potentially
 $K_6-C_5$-graphic.
 If $\pi_1 \in S$, then
 $\pi_1=$$(2)$,$(2^2)$,$(3^2)$,$(3^2,2)$.
 Therefore, $\pi=$$(6,3^6)$,$(7,3^7)$,$(5,3^7)$,$(6,3^8)$.
  It is easy to check that
 $\pi=$$(6,3^6)$,$(7,3^7)$ and $(6,3^8)$ are not potentially $K_6-C_5$-graphic
 and $\pi=(5,3^7)$, a contradiction to $d_1-5\neq 0$.\par

 If $1\leq d_2-3\leq (n-6)-(d_1-5)$, then let
 $\pi_1=(3^{(n-6)-(d_1-5)-(d_2-3)},\\2^{(d_1-5)+(d_2-3)})$.
 If $\pi_1 \notin S$, according to
 Lemma 2.6, then $\pi_1$ is graphic, hence $\pi$ is potentially
 $K_6-C_5$-graphic.If $\pi_1 \in S$,then
 $\pi_1 =$$(2)$,$(2^2)$,$(3^2)$,$(3^2,2)$.
Therefore, $\pi=(5,4,3^5)$,$(6,4,3^6)$,$(5^2,3^6)$,$(5,4,3^7)$. It
is easy to check that
 $\pi=(6,4,3^6)$ is not potentially $K_6-C_5$-graphic
 and $\pi=(5,4,3^5)$ ,$(5^2,3^6)$,
 $(5,4,3^7)$, a contradiction to $d_1-5\neq 0$.\par

 If \ $n-6 \geq d_2-3 > (n-6)-(d_1-5)$, \ then \ let
 $\pi_1=$ \ $(2^{(n-6)-(d_1-5)},\\1^{(d_2-3)+(d_1-5)-(n-6)},2^{(d_1-5)-(d_2-3)-(d_1-5)+(n-6)})$
 $=(2^{2(n-6)-(d_2-3)-(d_1-5)},\\1^{(d_1-5)+(d_2-3)-(n-6)})$.
If $\pi_1 \notin S$, according to Lemma 2.6, then $\pi_1$  is
graphic, hence  $\pi$ is potentially
 $K_6-C_5$-graphic.
 If $\pi_1 \in S$, then
 $\pi_1=$$(2)$,$(2^2)$. Therefore, $\pi=(5,4,3^5)$,$(6,4,3^6)$,$(5^2,3^6)$, a contradiction.
\par

 If $ d_2-3 > (n-6)$ and $d_2\leq n-1$, then $d_2=n-1$ or  $d_2=n-2$.
 If $d_2=n-1$, then $\pi=((n-1)^2,3^{n-2})$$(n\,\, is\,\, even)$ .
 It is easy to check that
 $\pi=((n-1)^2,3^{n-2})$$(n\,\, is\,\, even)$ is not potentially $K_6-C_5$-graphic.
 If $d_2=n-2$, then $\pi=(n-1,n-2,3^{n-2})$$(n \,\,is \,\,odd)$
 or $\pi=(n-2,n-2,3^{n-2})$$(n \,\,is \,\,even)$.
 It is easy to check that
 these are not potentially $K_6-C_5$-graphic.\par

   {\bf Lemma 2.11 }\ \ If $\pi=(d_1,3^5,2^{n-6})$ where $d_1\geq 5$,
 $n\geq 6$  and $\sigma(\pi)$ is even, then\par
 $(1)$ $\pi\in GS_n$;\par
 $(2)$  $\pi$ is potentially $K_6-C_5$-graphic if and only if
   $\pi\neq$$(5,3^5,2)$,$(5,3^5,2^2)$.

{\bf Proof:} $(1)$ is obvious\par
  $(2)$ If  $\pi=(d_1,3^5,2^{n-6})$ where $d_1\geq 5$,
 $n\geq 6$  and $\sigma(\pi)$ is even,
 then $\pi$ is potentially $K_6-C_5$-graphic if and only if
$\pi^*=(d_1-5,2^{n-6})$ is graphic. We consider the following
cases. \par

\textbf{Case 1:}  $d_1-5=0$\par
  If $\pi=(5,3^5,2^{n-6})$, then $\pi^*=(2^{n-6})$.
  If $\pi^*\notin S$, according to
 Lemma 2.6, then $\pi^*$  is graphic, hence  $\pi$ is potentially
 $K_6-C_5$-graphic.
 If $\pi^*\in S$, then
 $\pi^*=$$(2)$,$(2^2)$. Therefore, $\pi=(5,3^5,2)$,$(5,3^5,2^2)$.
 It is easy to check that
  these are not potentially $K_6-C_5$-graphic.\par

\textbf{Case 2:}  $1 \leq d_1-5\leq n-6$\par
  If $\pi=(d_1,3^5,2^{n-6})$, then let $\pi_1 =(2^{(n-6)-(d_1-5)},1^{d_1-5})$.
  Since $\pi_1 \notin S$, according to
 Lemma 2.6, $\pi_1$ is graphic, hence  $\pi$ is
potentially
 $K_6-C_5$-graphic.\par

   \par

   \par
   \section{ Main Theorems} \par
   \par
 \ \ \ \  \textbf{\noindent Theorem 3.1}\ \  Let
 $\pi=(d_1,d_2,\cdots,d_n)$ $\in GS_n$ and
 $n\geq6$. Then $\pi$ is potentially $K_6-C_5$-graphic if and only
if the following conditions hold:
\par
  $(1)$ $d_1\geq5,d_6\geq3$;\par
  $(2)$ If $\pi=(d_1,d_2,d_3,3^i,2^j,1^{n-i-j-3})$ where $i\geq3$ and $d_3\geq5$, then
        $d_1+d_2+d_3\leq n+2i+j+1$;\par
  $(3)$ If $\pi=(d_1,d_2,3^4,2^j,1^{n-j-6})$ where $j\geq0$ and $d_2\geq5$, then
        $d_1+d_2\leq n+j+2$;\par
  $(4)$ If $\pi=(d_1,d_2,d_3,3^i,2^j,1^{n-i-j-3})$ where $n\geq 8$, $i\geq4$
        and $d_3\geq4$, then\par
        $(i)$ if $n=i+3$, $d_1=n-1$ and $d_3\geq 5$, then $d_2\leq n-2$;\par
        $(ii)$ if $n\geq i+j+4$, $d_1\geq i+j+3$, $d_2\geq max\{d_3+2, i+j+2\}$ and $d_3=4$, then $d_1+d_2\leq n+i+j$;\par
        $(iii)$ if $n\geq i+j+4$, $d_1=n-1$, $d_2\geq d_3+2$ and $d_3\geq 5$,  then $d_1+d_2\leq n+i+j$;\par
        $(iv)$ if $n= i+j+3$, $j\geq 1$, $d_1= n-1$ and $d_2\geq d_3+2$, then $d_2\leq n-2$;\par
  $(5)$ If $\pi=(d_1,d_2,3^i,2^j,1^{n-i-j-2})$ where $i\geq5$ and
        $d_2\geq5$ , then\par
        $(i)$ if $n\geq i+1$, $j\geq 2$ or $j=0$, $d_1\geq i+j$, then $d_1+d_2\leq n+i+j-2$;\par
        $(ii)$ if $d_1= i+j+1$ and $d_1\geq n-2$, then $d_2\leq n-3$;\par
        $(iii)$ if $d_1= i+j=n-2$ and $i\geq 6$, then $d_2\leq n-3$;\par
  $(6)$ $\pi\neq((n-1)^2,4,3^{n-3}) (n\geq 7\ \  and \ \ \ n\ \  is \ \ odd) $;\par
  $(7)$ $\pi\neq(n-1,3^6,1^{n-7}) (n\geq 7)$,$(n-1,3^7,1^{n-8}) (n\geq 8)$,$(5,4,3^5)$,$(5^2,3^5,1)$,\\
  $(5,4^2,3^5)$, $(5,4,3^7)$, $(6,3^7,1)$, $(5^2,3^6)$, $(5^2,4,3^4)$, $(5,3^7)$, $(6,4,3^6)$, $(5^2,3^4,2)$,
  $(5,3^6,1)$, $(5,3^6,2,1)$, $(5,3^7,2)$, $(5,3^8,1)$, $(5,4,3^5,2)$, $(5,4,3^6,1)$, $(5,3^5,2)$, $(5,3^5,2^2)$,
  $(6,5,3^5)$, $(5,3^7,1^2)$, $(6,3^8)$, $(6^2,3^4,2)$.\par
\par
{\bf Proof:} First we show the conditions $(1)$ to $(7)$ are
necessary conditions for $\pi$ to be potentially
$K_6-C_5$-graphic. Assume that $\pi$ is potentially
$K_6-C_5$-graphic.\par
  $(1)$ is obvious.\par

  $(2)$ If $\pi=(d_1,d_2,d_3,3^i,2^j,1^{n-i-j-3})$ where $i\geq3$ and $d_3\geq5$, then
  $(d_1+d_2+d_3-11)-2\leq 3\times (i-3)+2\times j+1\times (n-i-j-3)$,
  i.e., $d_1+d_2+d_3\leq n+2i+j+1$.\par

  $(3)$ If $\pi=(d_1,d_2,3^4,2^j,1^{n-j-6})$ where $j\geq0$ and $d_2\geq5$, then
   $(d_1+d_2-8)\leq 2\times j+(n-j-6)$, i.e., $d_1+d_2\leq n+j+2$.\par

  $(4)$   If $\pi=(d_1,d_2,d_3,3^i,2^j,1^{n-i-j-3})$ where  $n\geq 8$, $i\geq4$
        and $d_3\geq4$, then\par
        $(i)$ if $n=i+3$, $d_1=n-1$ and $d_3\geq 5$, then  $d_2-3-1\leq n-6$ , i.e., $d_2\leq n-2$;\par

        $(ii) $ if $n\geq i+j+4$, $d_1\geq i+j+3$, $d_2\geq max\{d_3+2, i+j+2\}$ and $d_3=4$, then
        $(d_2-3)-1\leq (i-3)+j+(n-i-j-3)+(i-3)+j-(d_1-5)$
        or $(d_1-3)-1\leq (i-3)+j+(n-i-j-3)+(i-3)+j-(d_2-5)$,  i.e., $d_1+d_2\leq n+i+j$.\par

        $(iii)$ if $n\geq i+j+4$, $d_1=n-1$, $d_2\geq d_3+2$ and $d_3\geq 5$, then
        $(d_2-3)-1\leq (i-3)+j+(n-i-j-3)+(i-3)+j-(d_1-5)$, i.e., $d_1+d_2\leq n+i+j$
        or $(d_1-3)-1\leq (i-3)+j+(n-i-j-3)+(i-3)+j-(d_2-5)$, i.e.,  $d_1+d_2\leq n+i+j$;\par

        $(iv)$ if $n= i+j+3$, $j\geq 1$, $d_1= n-1$, $d_2\geq d_3+2$, then
        $d_2-3-1\leq (i-3)+j+(n-i-j-3) = n-6$ , i.e., $d_2\leq n-2$.\par

  $(5)$ If $\pi=(d_1,d_2,3^i,2^j,1^{n-i-j-2})$ where $i\geq5$ and $d_2\geq5$, then\par
        $(i)$ if $n\geq i+1$, $j\geq 2$ or $j=0$, $d_1\geq i+j$,  then
        $d_2-3\leq (i-4)+j+(n-i-j-2)-[(d_1-5)-(i-4)-j]$
        or $d_1-5\leq (i-4)+j+(n-i-j-2)-[(d_2-3)-(i-4)-j]$, i.e., $d_1+d_2\leq n+i+j-2$;\par

        $(ii)$ if $d_1= i+j+1$ and $d_1\geq n-2$, then $\pi^*=(d_1-5,d_2-3,3^{i-4},2^j,1^{n-i-j-2})$.
        Therefore, $d_2-3\leq (i-4)+j+(n-i-j-2)$, i.e., $d_2\leq n-3$.\par

        $(iii)$ if $d_1= i+j=n-2$ and $i\geq 6$, then $\pi^*=(d_1-5,d_2-3,3^{i-4},2^{n-i-2})$.
        Therefore, $d_2-3\leq (i-4)+(n-i-2)$, i.e. $d_2\leq n-3$;\par

  $(6)$ If $\pi=((n-1)^2,4,3^{n-3}) (n\geq 7 \ \ and \ \ n \ \ is\ \  odd) $ is potentially
  $K_6-C_5$-graphic,
        then  $\pi^*=(n-6,n-4,1,3^{n-6})$ is graphic, hence $n-4-1\leq n-6$, a contradiction.\par

  $(7)$ If $\pi=(n-1,3^6,1^{n-7}) (n\geq 7)$ is potentially
  $K_6-C_5$-graphic, then  $\pi^*=(n-6,3,1^{n-7})$ is graphic, hence
  $\pi_1=(2)$$\in$$GS_n$, a contradiction.\par

If $\pi=(n-1,3^7,1^{n-8}) (n\geq 8) $is potentially
  $K_6-C_5$-graphic, then  $\pi^*=(n-6,3^2,1^{n-8})$ is graphic, hence
  $\pi_1=(2^2)$$\in$$GS_n$, a contradiction.\par
If $\pi=$$(5,4,3^5)$,$(5^2,3^5,1)$,$(5,4^2,3^5)$,$(5,4,3^7)$,$(6,3^7,1)$,$(5^2,3^6)$,$(5^2,4,3^4)$,\\
$(5,3^7)$,$(6,4,3^6)$,$(5^2,3^4,2)$,$(5,3^6,1)$,$(5,3^6,2,1)$,$(5,3^7,2)$,$(5,3^8,1)$,$(5,4,3^5,2)$,\\
$(5,4,3^6,1)$,$(5,3^5,2)$,$(5,3^5,2^2)$,$(6,5,3^5)$,$(5,3^7,1^2)$,$(6,3^8)$,
$(6^2,3^4,2)$,
          then it is easy to see that they are not potentially $K_6-C_5$-graphic .\par

  Now we prove the sufficient condition. Suppose the graphic sequence
  $\pi$ satisfies the conditions (1) to (7). Our proof is by induction
  on $n$.\par
  We first prove the base case where $n=6$.
  Then $\pi$ is one of the following:
  $(5^6)$,$(5^4,4^2)$ ,$(5^3,4^2,3)$,$(5^3,3^3)$,$(5^2,4^4)$,$(5^2,4^2,3^2)$, $(5^2,3^4)$,$(5,4^4,3)$,
  $(5,4^2,\\3^3)$, $(5,3^5)$ .
  It is easy to check that others  are potentially $K_6-C_5$-graphic except for $(5^3,3^3)$ and $(5^2,3^4)$.
  Now suppose that the sufficiency holds for $n-1(n\geq7)$,
  we will show that $\pi$ is potentially  $K_6-C_5$-graphic in terms of the following cases:\par

\textbf{Case 1:} $d_n\geq 6$\par Consider
$\pi^\prime=(d_1^\prime,d_2^\prime,\cdots,d_{n-1}^\prime)$ where
 $d^\prime_{n-1}\geq 5$.\par
 It is easy to check that
$\pi^\prime$ satisfies (1) to (7). Then by the induction
hypothesis, $\pi^\prime$ is potentially $K_6-C_5$-graphic, and
hence so is $\pi$.
\par

\textbf{Case 2:} $d_n=5$ \par Consider
$\pi^\prime=(d_1^\prime,d_2^\prime,\cdots,d_{n-1}^\prime)$ where
$d^\prime_1\geq 5$ and $d^\prime_{n-1}\geq 4$.\par
 It is easy to
check that $\pi^\prime$ satisfies (1) to (7).
   Then by the induction hypothesis,
$\pi^\prime$ is potentially $K_6-C_5$-graphic, and hence so is
$\pi$.

\textbf{Case 3:} $d_n=4$\par Consider
$\pi^\prime=(d_1^\prime,d_2^\prime,\cdots,d_{n-1}^\prime)$ where
$d_1^\prime\geq4$ and $d_{n-1}^\prime\geq 3$.\par

 It is easy to check that $\pi^\prime$ satisfies $(2)$ to $(7)$.
 If $\pi^\prime$ satisfies
$(1)$, then by the induction hypothesis, $\pi^\prime$ is potentially
$K_6-C_5$-graphic, and hence so is $\pi$.\par

If $\pi^\prime$ does not satisfy $(1)$, then
 $d_1^\prime=4$. Therefore, $d_5=4=d_6=\cdots=d_{n-1}$.
Then $\pi=$$(5^2,4^{n-2})$$(n\geq 7)$ or
$\pi=$$(5^4,4^{n-4})$$(n\geq 7)$. According to Lemma 2.8 and 2.7,
$\pi=$$(5^2,4^{n-2})$$(n\geq 7)$ and $\pi=$$(5^4,4^{n-4})$$(n\geq
7)$ holding the conditions of Theorem 3.1 are potentially
$K_6-C_5$-graphic.\par

\textbf{Case 4:} $d_n=3$\par Consider
$\pi^\prime=(d_1^\prime,d_2^\prime,\cdots,d_{n-1}^\prime)$ where
$d_1^\prime\geq 4$, $d_{n-3}^\prime\geq3$ and $d_{n-1}^\prime\geq
2$.\par

It is easy to check that $\pi^\prime$ satisfies $(2)$.
 If $\pi^\prime$ satisfies
$(1)$,$(3)$ to $(7)$, then by the induction hypothesis,
$\pi^\prime$ is potentially $K_6-C_5$-graphic, and hence so is
$\pi$.\par

If $\pi^\prime$ does not satisfy $(1)$, then
 $d^\prime_1=4$ or $d^\prime_6\leq 2$.
If $d^\prime_1=4$, then $d_1=5$ and $3 \leq d_4 \leq 4$. Therefore,
$\pi=$$(5^i,4^j,3^{n-i-j})$$(1 \leq i \leq 3,j \geq 0,n \geq 7)$.

According to Lemma 2.7 to 2.9, $\pi=$$(5^i,4^j,3^{n-i-j})$$(1 \leq
i\leq 3,j \geq 0,n \geq 7)$ holding the conditions of Theorem 3.1
is potentially $K_5-C_4$-graphic.\par

If $d^\prime_6=2$, then $d_3=3=\cdots=d_n$ .Therefore,
$\pi=$$(d_1,d_2,3^{n-2})$$(d_1\geq 5, d_2\geq 3, n\geq 7)$.
According to Lemma 2.10,
$\pi=$$(d_1,d_2,3^{n-2})$$(d_1\geq5,d_2\geq3)$ holding the
conditions of Theorem 3.1 is potentially $K_6-C_5$-graphic.\par

If $\pi^\prime$ does not satisfy $(3)$, then $n=7$ and
$\pi^\prime=(5^2,3^4)$, hence $\pi=$$(6,6,4,\\3^4)$, a
contradiction to the condition $(6)$.\par

If $\pi^\prime$ does not satisfy $(4)$, then $\pi^\prime$
satisfies $(4)$$(ii)$,  $(4)$$(iii)$ and  $(4)$$(iv)$.\par

If $\pi^\prime$ does not satisfy $(4)$$(i)$, then
$\pi^\prime=((n-1)^2,d^\prime_3,3^{n-3})$ where $i \geq4$ and
$d^\prime_3\geq 5$, hence $\pi=(n^2,d^\prime_3+1,3^{n-2})$, a
contradiction to the condition $(4)(i)$.\par

If $\pi^\prime$ does not satisfy $(5)$, then $\pi^\prime$ satisfies
$(5)$$(i)$.\par

If $\pi^\prime$ does not satisfy $(5)$$(ii)$, then
$\pi^\prime=(n-1,n-2,3^{n-j-2},2^j)$ or
$\pi^\prime=((n-1)^2,3^{n-j-2},2^j)$ where $n-j-2\geq 5$ and $j\leq
1$. If $\pi^\prime=(n-1,n-2,3^{n-2})$, then
$\pi=(n,n-1,4,3^{n-2})$($n\geq 7$ and n is odd) is potentially
$K_6-C_5$-graphic . If $\pi^\prime=(n-1,n-2,3^{n-3},2)$, then
$\pi=(n,n-1,3^{n-1})$, a contradiction to the condition $(5)(i)$. If
$\pi^\prime=((n-1)^2,3^{n-2})$, then $\pi=(n^2,4,3^{n-2})$, a
contradiction to the condition $(6)$. If
$\pi^\prime=((n-1)^2,3^{n-3},2)$, then $\pi=(n^2,3^{n-1})$, a
contradiction to the condition $(5)(ii)$.\par

If $\pi^\prime$ does not satisfy $(5)$$(iii)$, then
$\pi^\prime=((n-2)^2,3^{n-2})$($n\geq 8$) or
$\pi^\prime=((n-2)^2,3^{n-3},2)$($n\geq 9$). Therefore,
$\pi=((n-1)^2,4,3^{n-2})$($n\geq 8$) is potentially
$K_6-C_5$-graphic  or $\pi=((n-1)^2,3^{n-1})$($n\geq 9$) which
contradicts  the condition $(5)(iii)$.\par

If $\pi^\prime$ does not satisfy $(6)$, then
$\pi=$$(n^2,5,3^{n-2})$ which contradicts  the condition $(4)(i)$
or $\pi=$$(n^2,4^2,3^{n-3})$ is potentially $K_6-C_5$-graphic
.\par

If $\pi^\prime$ does not satisfy $(7)$, then
  $\pi^\prime=(6,3^6)$,$(7,3^7)$,$(5,4,3^5)$,$(5,4^2,3^5)$, $(5,4,3^7)$,
   $(5^2,3^6)$,$(5^2,4,3^4)$,$(5,3^7)$,$(6,4,3^6)$,$(5^2,3^4,2)$,$(5,3^7,2)$,$(5,4,3^5,2)$,\\
  $(5,3^5,2)$,\,\, $(5,3^5,2^2)$,\,\, $(6,5,3^5)$,\,\, $(6,3^8)$,\,\,$(6^2,3^4,2)$, hence\,\, $\pi=$$(7,4^2,3^5)$,\,\, $(8,4^2,3^6)$,
  $(6,5,4,3^5)$,$(6,4^3,3^4)$,$(6,5^2,3^6)$,$(6,5,4^2,3^5)$,$(6,4^4,3^4)$,$(6,5,4,3^7)$,\\$(6,4^3,3^6)$,
  $(6^2,4,3^6)$, $(6^2,5,3^5)$, $(6^2,4^2,3^4)$, $(6,4^2,3^6)$, $(7,5,4,3^6)$, $(7,4^3,3^5)$, $(6^2,3^6)$,
  $(6,4,3^8)$, $(6,5,3^7)$, $(6,4,3^6)$,$(6,3^8)$,$(7,6,4,3^5)$,$(7,4^2,3^7)$,$(7^2,4,3^5)$,\\$(7^2,3^6)$.
  It is easy to check that
 others are  potentially $K_6-C_5$-graphic except
 $(6^2,3^6)$, $(6,4,3^6)$, $(6,3^8)$, $(7^2,4,3^5)$ and
$(7^2,3^6)$. Then $\pi=$$(6^2,3^6)$,$(6,4,\\3^6)$
$(6,3^8)$,$(7^2,4,3^5)$ and $(7^2,3^6)$, a contradiction to the
conditions $(5)(iii)$, $(7)$,  $(4)(iv)$, $(6)$ and $(5)(ii)$.
\par

\textbf{Case 5:} $d_n=2$\par Consider
$\pi^\prime=(d_1^\prime,d_2^\prime,\cdots,d_{n-1}^\prime)$ where
$d_1^\prime\geq4$, $d_5^\prime\geq 3$ and
$d_{n-1}^\prime\geq2$.\par

It is easy to check that $\pi^\prime$ satisfies $(2)$ and $(3)$.
 If $\pi^\prime$ satisfies
$(1)$ and $(4)$ to $(7)$, then by the induction hypothesis,
$\pi^\prime$ is potentially $K_6-C_5$-graphic, hence so is
$\pi$.\par

If $\pi^\prime$ does not satisfy $(1)$, then
 $d^\prime_1=4$ or $d^\prime_6\leq2$.
If $d^\prime_1=4$, then $d_1=5$ and $d_3\leq 4$. Therefore,
$\pi=$$(5^2,4^i,3^j,2^{n-2-i-j})$$(i+j\geq4,n\geq 7)$ or
$\pi=$$(5,4^i,3^j,2^{n-1-i-j})$$(i+j\geq5,n\geq 7)$.

According to Lemma 2.8 and 2.9,
 $\pi=$$(5^2,4^i,3^j,2^{n-2-i-j})$$(i+j\geq4,n\geq 7)$ and
$\pi=$$(5,4^i,3^j,2^{n-1-i-j})$$(i+j\geq5,n\geq 7)$ holding the
conditions of Theorem 3.1 are potentially $K_6-C_5$-graphic.\par

If $d^\prime_6=2$, then $d_2=3=d_3=d_4=d_5=d_6$ and
$d_7=\cdots=d_n$. Therefore,
$\pi=$$(d_1,3^5,2^{n-6})$$(d_1\geq5,n\geq7)$ .

According to Lemma 2.11,
 $\pi=$$(d_1,3^5,2^{n-6})$$(d_1\geq5,n\geq7)$ holding the conditions of Theorem
3.1 is potentially $K_6-C_5$-graphic.\par

If $\pi^\prime$ does not satisfy $(4)$, then $\pi^\prime$
satisfies $(4)$$(ii)$ and $(4)$$(iii)$.\par

If $\pi^\prime$ does not satisfy $(4)$$(i)$, then
$\pi^\prime=((n-1)^2,d^\prime_3,3^{n-3})$ where $i \geq4$ and
$d^\prime_3\geq 5$. If $d^\prime_3\neq n-2$, then
$\pi=(n^2,d^\prime_3,3^{n-3},2)$, a contradiction to the condition
(4)(i). If $d^\prime_3= n-2$, then $\pi=(n^2,n-2,3^{n-3},2)$ which
contradicts the condition (2) or $\pi=(n,(n-1)^2,3^{n-3},2)$ which
also contradicts the condition (2).\par

If $\pi^\prime$ does not satisfy $(4)$$(iv)$, then
$\pi^\prime=((n-1)^2,d^\prime_3,3^i,2^{n-i-3})$ where $i \geq 5$ ,
$d^\prime_3\geq4$ and $d^\prime_2 \geq d^\prime_3+2$.  Therefore,
$\pi=(n^2,d^\prime_3,3^i,2^{n-i-2})$, a contradiction to the
condition $(4)(iv)$.

If $\pi^\prime$ does not satisfy $(5)$, then $\pi^\prime$ satisfies
$(5)$$(i)$.\par

If $\pi^\prime$ does not satisfy $(5)$$(ii)$ where $i \geq4$ and
$d^\prime_2\geq 5$, then
$\pi^\prime=(d^\prime_1,d^\prime_2,3^i,2^{n-i-2})$. Therefore,
$\pi=(d^\prime_1+1,d^\prime_2+1,3^i,2^{n-i-1})$, a contradiction
to the condition $(5)(ii)$.\par

If $\pi^\prime$ does not satisfy $(5)$$(iii)$, then
$\pi^\prime=((n-2)^2,3^i,2^{n-i-2})$($i\geq 6$). Therefore,
$\pi=((n-1)^2,3^i,2^{n-i-2+1})$($i\geq 6$), a contradiction to the
condition $(5)(iii)$.\par

If $\pi^\prime$ does not satisfy $(6)$, then
$\pi=$$(n^2,4,3^{n-3},2)$, a contradiction to the condition
$(4)(iii)$.\par

If $\pi^\prime$ does not satisfy $(7)$, then
  $\pi^\prime=(6,3^6)$,$(7,3^7)$,$(5,4,3^5)$,
 $(5,4^2,3^5)$, $(5,4,3^7)$,
   $(5^2,3^6)$,$(5^2,4,3^4)$,$(5,3^7)$,
  $(6,4,3^6)$,$(5^2,3^4,2)$,
  $(5,3^7,2)$,$(5,4,3^5,2)$,\\
  $(5,3^5,2)$,\, $(5,3^5,2^2)$,\, $(6,5,3^5)$,\, $(6,3^8)$, $(6^2,3^4,2)$,
 hence $\pi=$$(7,4,3^5,2)$,\, $(8,4,3^6,2)$, $(6,4^2,3^4,2)$,
 $(6,5,3^5,2)$, $(6,5,4,3^5,2)$, $(5^3,3^5,2)$, $(6,4^3,3^4,2)$, $(6,5,3^7,2)$, $(6,4^2,3^6,2)$,
 $(6^2,3^6,2)$, $(6^2,4,3^4,2)$, $(6,5^2,3^4,2)$, $(6,4,3^6,2)$, $(7,5,3^6,2)$, $(7,4^2,3^5,2)$,
 $(6^2,3^4,2^2)$, $(6,4,3^6,2^2)$, $(6,5,3^5,2^2)$, $(6,4^2,3^4,2^2)$, $(6,4,3^4,2^2)$, $(6,4,3^4,2^3)$,
 $(7,6,3^5,2)$,$(7,4,3^7,2)$,$(7^2,3^4,2^2)$.
It is easy to check that others are  potentially $K_6-C_5$-graphic
except  $(6^2,3^4,2^2)$, $(7,6,3^5,\\2)$ and $(7^2,3^4,2^2)$. Then
$\pi=$$(6^2,3^4,2^2)$,$(7,6,3^5,2)$ and $(7^2,3^4,2^2)$, a
contradiction to the conditions $(5)(iii)$ and $(5)(ii)$.\par

\textbf{Case 6:} $d_n=1$ \par
Consider
$\pi^\prime=(d_1^\prime,d_2^\prime,\cdots,d_{n-1}^\prime)$ where
$d_1^\prime\geq4$, $d_6^\prime\geq3$ and
$d_{n-1}^\prime\geq1$.\par

It is easy to check that $\pi^\prime$ satisfies $(2)$ and $(3)$.
 If $\pi^\prime$ satisfies
$(1)$ and $(4)$ to $(7)$, then by the induction hypothesis,
$\pi^\prime$ is potentially $K_6-C_5$-graphic, and hence so is
$\pi$.\par

If $\pi^\prime$ does not satisfy $(1)$, then
 $d^\prime_1=4$.
If $d^\prime_1=4$, then $d_1=5$. Therefore,
$\pi=$$(5,4^i,3^j,2^k,1^{n-i-j-k-1})$$(i+j\geq5,n\geq 7)$ .

According to Lemma 2.9,
$\pi=$$(5,4^i,3^j,2^k,1^{n-i-j-k-1})$$(i+j\geq5,n\geq 7)$ holding
the conditions of Theorem 3.1 is potentially
$K_6-C_5$-graphic.\par

If $\pi^\prime$ does not satisfy $(4)$$(i)$, then
$\pi^\prime=((n-1)^2,d^\prime_3,3^{n-3})$ where $i \geq4$ and
$d^\prime_3\geq 5$.\par

If $d^\prime_3\neq n-2$, then $\pi=(n,n-1,d^\prime_3,3^{n-3},1)$,
a contradiction to the condition $(4)(iii)$. If $d^\prime_3= n-2$,
then $\pi=(n,n-1,n-2,3^{n-3},1)$
 or $\pi=((n-1)^3,3^{n-3},1)$ which also contradicts
the conditions $(4)$$(iii)$ and $(2)$.\par

If $\pi^\prime$ does not satisfy $(4)$$(ii)$, then
$\pi^\prime=(d^\prime_1,d^\prime_2,d^\prime_3,3^i,2^j,1^{n-i-j-3})$
where $i \geq4$ , $d^\prime_3\geq4$ and $d^\prime_2 \geq
d^\prime_3+2$.\par

If $d^\prime_1 \neq d^\prime_2+1$,
 then $\pi=(d^\prime_1+1,d^\prime_2,d^\prime_3,3^i,2^j,1^{n-i-j-3})$, a contradiction
to the condition $(4)(ii)$. If $d^\prime_1 = d^\prime_2+1$,
 then
 $\pi=(d^\prime_1+1,d^\prime_2,d^\prime_3,3^i,2^j,1^{n-i-j-3})$ which contradicts  the condition $(4)(ii)$
 or $\pi=(d^\prime_1,d^\prime_2+1,d^\prime_3,3^i,2^j,1^{n-i-j-3})$ which also contradicts  the condition $(4)(ii)$.\par

If $\pi^\prime$ does not satisfy $(4)$$(iii)$, then
$\pi^\prime=(d^\prime_1,d^\prime_2,d^\prime_3,3^i,2^j,1^{n-i-j-3})$
and $i \geq4$ , $d^\prime_3\geq 5$, $d^\prime_2 \geq
d^\prime_3+2$.\par

If $d^\prime_1 \neq d^\prime_2+1$, then
$\pi=(d^\prime_1+1,d^\prime_2,d^\prime_3,3^i,2^j,1^{n-i-j-3})$,
which contradicts  the condition $(4)(iii)$. If $d^\prime_1 =
d^\prime_2+1$, then
$\pi=(d^\prime_1+1,d^\prime_2,d^\prime_3,3^i,2^j,\\1^{n-i-j-3})$
or $\pi=(d^\prime_1,d^\prime_2+1,d^\prime_3,3^i,2^j,1^{n-i-j-3})$,
which contradicts  the condition $(4)(iii)$.\par

If $\pi^\prime$ does not satisfy $(4)$$(iv)$, then
$\pi^\prime=((n-1)^2,d^\prime_3,3^i,2^j,1^{n-i-j-3})$ and
$i\geq4$, $d^\prime_2 \geq d^\prime_3+2$. Then
$\pi=(n,n-1,d^\prime_3,3^i,2^j,1^{n-i-j-2})$, which contradicts
the conditions $(4)(ii)$ and $(4)(iii)$.\par

If $\pi^\prime$ does not satisfy $(5)$$(i)$, then
$\pi^\prime=(d^\prime_1,d^\prime_2,3^i,2^j,1^{n-i-j-2})$ where
$i\geq 5$ and $d_2^\prime \geq 5$. If $d^\prime_1\neq
d^\prime_2-1$,
 then $\pi=(d^\prime_1+1,d^\prime_2,3^i,2^j,1^{n-i-j-1})$, a
contradiction to the condition $(5)(i)$. If $d^\prime_1=
d^\prime_2-1$, then
$\pi=(d^\prime_1+1,d^\prime_2,3^i,2^j,1^{n-i-j-1})$ which
contradicts  the condition $(5)(i)$ or
$\pi=(d^\prime_1,d^\prime_2+1,3^i,2^j,1^{n-i-j-1})$ which also
contradicts  the condition $(5)(i)$ .\par

If $\pi^\prime$ does not satisfy $(5)$$(ii)$, then
$\pi^\prime=(n-1,n-2,3^i,2^{n-i-2})$,$((n-1)^2,3^i,2^{n-i-2})$ or
$((n-2)^2,3^i,2^{n-i-3},1)$. If
$\pi^\prime=(n-1,n-2,3^i,2^{n-i-2})$, then
$\pi=(n,n-2,3^i,2^{n-i-2},1)$ which contradicts  the condition
$(5)(i)$ or $\pi=((n-1)^2,3^i,2^{n-i-2},1)$ which contradicts the
condition $(5)(ii)$. If $\pi^\prime=((n-1)^2,3^i,2^{n-i-2})$, then
$\pi=(n,n-1,3^i,2^{n-i-2},1)$, a contradiction to the condition
$(5)(i)$. If $\pi^\prime=((n-2)^2,3^i,2^{n-i-3},1)$, then
$\pi=(n-1,n-2,3^i,2^{n-i-2},1^2)$, a contradiction to the
condition $(5)(i)$.\par

If $\pi^\prime$ does not satisfy $(5)$$(iii)$, then
$\pi^\prime=((n-2)^2,3^i,2^{n-i-2})$($i\geq 6$). Therefore,
$\pi=(n-1,n-2,3^i,2^{n-i-2},1)$ ($i\geq 6$) is potentially
$K_6-C_5$-graphic.\par

If $\pi^\prime$ does not satisfy $(6)$, then
$\pi=$$(n,n-1,4,3^{n-3},1)$, a contradiction to the condition
$(4)(ii)$.\par

If $\pi^\prime$ does not satisfy $(7)$, then
  $\pi^\prime=$$(n-1,3^6,1^{n-7}) (n\geq 7) $,$(n-1,3^7,\\1^{n-8}) (n\geq 8) $, $(5,4,3^5)$,
  $(5^2,3^5,1)$, $(5,4^2,3^5)$, $(5,4,3^7)$, $(6,3^7,1)$, $(5^2,3^6)$, $(5^2,4,3^4)$, $(5,3^7)$,$(6,4,3^6)$,
  $(5^2,3^4,2)$,$(5,3^6,1)$,$(5,3^6,2,1)$,$(5,3^7,2)$,$(5,3^8,1)$, $(5,4,3^5,2)$, $(5,4,3^6,1)$,
  $(5,3^5,2)$, $(5,3^5,2^2)$, $(6,5,3^5)$, $(5,3^7,1^2)$,
  $(6,3^8)$, $(6^2,3^4,2)$.
Hence $\pi=$$(n,3^6,1^{n-6}) (n\geq 7)$,$(n,3^7,1^{n-7})(n\geq8)$,
$(6,4,3^5,1)$, $(5^2,3^5,1)$, $(6,5,3^5,1^2)$, $(6,4^2,3^5,1)$,
$(5^2,4,3^5,1)$, $(6,4,3^7,1)$, $(5^2,3^7,1)$,
  $(7,\\3^7,1^2)$, $(6,5,3^6,1)$, $(6,5,4,3^4,1)$, $(5^3,3^4,1)$, $(6,3^7,1)$, $(7,4,3^6,1)$,
  $(6,5,3^4,2,\\1)$, $(6,3^6,1^2)$, $(6,3^6,2,1^2)$, $(6,3^7,2,1)$, $(6,3^8,1^2)$, $(6,4,3^5,2,1)$,
  $(5^2,3^5,2,1)$, $(6,4,3^6,1^2)$, $(5^2,3^6,1^2)$, $(6,3^5,2,1)$,\, $(6,3^5,2^2,1)$, $(7,5,3^5,1)$,
  $(6^2,3^5,1)$, $(6,3^7,1^3)$, $(7,3^8,1)$,$(7,6,3^4,2,1)$. It is easy to check that
 others are  potentially $K_6-C_5$-graphic except
  $\pi=$$(n,3^6,1^{n-6}) (n\geq 7)$, $(n,3^7,1^{n-7})(n\geq8)$, $(5^2,3^5,1)$, $(6,3^7,1)$,
  $(7,5,3^5,1)$,  $(6^2,3^5,1)$ and $(7,6,3^4,2,1)$. Then $\pi=$$(n,3^6,1^{n-6}) (n\geq 7)$, $(n,3^7,1^{n-7})(n\geq8)$,
$(5^2,3^5,1)$, $(6,3^7,1)$, $(7,5,3^5,1)$, $(6^2,3^5,1)$ and
$(7,6,3^4,2,1)$, a contradiction to the conditions $(7)$, $(5)$
and $(4)$.\par

\par
\vspace{0.5cm}
\par
\textbf{Corollary 3.2} For $n\geq6$, $\sigma(K_6-C_5,n)=6n-10$.
\par
{\bf Proof:}  First we claim $\sigma(K_6-C_5,n)\geq 6n-10$ for $n
\geq6$. We would like to show there exists $\pi_1$ with
$\sigma(\pi_1)=6n-12$ such that $\pi_1$ is not potentially
$K_6-C_5$-graphic. Let $\pi_1=((n-1)^3,3^{n-3})$. It is easy to see
that $\sigma(\pi_1)=6n-12$ and $\pi_1$ is not potentially
$K_6-C_5$-graphic by Theorem 3.1.\par

  Now we easily show if $\pi$ is an $n$-term $(n\geq6)$ graphic sequence
with $\sigma(\pi) \geq 6n-10$, then there exists a realization of
$\pi$ containing a $K_6-C_5$.\par

  $(1)$ We claim that $d_1\geq 5$ and $d_6 \geq 3$. By way of
  contradiction, we assume that $d_1\leq 4$ or $d_6 \leq 2$. If $d_1\leq
  4$, then $\sigma(\pi_1)$$\leq 4\times(n-1)$$< 6n-12(n\geq 6)$, a
  contradiction.  If $d_6\leq 2$, then $\sigma(\pi_1)$$=\sum_{i=1}^5 d_i$+$\sum_{i=6}^n d_i$$\leq $
  $5\times(5-1)$+$\sum_{i=6}^n  min(5,d_i)$+$\sum_{i=6}^n d_i$=20+2$\sum_{i=6}^n d_i$
  $\leq$ 4n $<6n-10(n\geq 6)$, a
  contradiction.\par
  $(2)$ If $\pi=(d_1,d_2,d_3,3^i,2^j,1^{n-i-j-3})$ where $i\geq3$ and
  $d_3\geq5$, then
  $\sigma(\pi)=$$d_1+d_2+d_3+3i+2j+n-i-j-3$ $\leq 3\times(n-1)+3\times(n-3)=6n-12$, a contradiction.
  Then $\pi$ with $\sigma(\pi) \geq 6n-10$ satisfies the condition $(2)$ of
 the Theorem 3.1.\par
  Similar to prove that $\pi$ with $\sigma(\pi) \geq 6n-10$ satisfies the conditions $(3)$ to $(7)$ of
 the Theorem 3.1. According to the Theorem 3.1, $\pi$ is  potentially
$K_6-C_5$-graphic.
\par

\section*{Acknowledgment}{}
 \ \ \ \ \  The authors wish to thank Professor Gang Chen, R.J.
Gould,  Jiongsheng Li, Rong Luo,  John R. Schmitt,  Zi-Xia Song,
Amitabha Tripathi, Jianhua Yin and Mengxiao Yin for sending some
their papers to us.

\end{document}